# Ridesharing Evacuation Model of Disaster Response

## Abstract ID: 797308

Lingyu Meng, Zhijie (Sasha) Dong
Texas State University
San Marcos, Texas, USA

## Abstract

Timely evacuation is crucial to disaster response, as people can avoid suffering and loss of lives when a major disaster happens. With the development of sharing economy, ridesharing has the advantage of reducing congestion, saving travel time, and optimizing transportation mode to improve disaster evacuation efficiency. The paper proposes to integrate the concept of "ridesharing" into evacuation and develops a mixed-integer programming model for this problem. A real-world case study based on Houston is used to validate the proposed model. A series of instances are designed to compare the evacuation efficiency using two indicators, evacuation percentage (EP) and average travel distance (ATD). Results reveal that increasing the number of vehicles to help carless individuals might not be the most efficient method in this model. Moreover, this model offers a specific response strategy based on different disaster scales, which not only develops a better evacuation plan for the people but also provides relief agencies insights on resource utilization.

## 1. Introduction

As one of the most significant activities of disaster response, evacuation can mitigate the adverse effects of natural disasters, such as suffering and loss of lives for individuals. However, resources like vehicles are usually insufficient. Therefore, improving evacuation efficiency with limited resource is of great importance.

The field of disaster evacuation has been extensively studied. Most existing studies have focused on public transport evacuation, such as proposing public transport evacuation plans for people without access to a private vehicle [1] and analyzing public transport plans in different disaster scenarios [2]. These studies assumed that people have already been gathered at specific places before the evacuation, and then they developed evacuation routes or allocated vehicles - however, few studies considered how to use private vehicles to transfer people for evacuation

Ridesharing is an emerging topic in which individual travelers share a private vehicle for a trip with others that have similar itineraries and time schedules. Most existing studies are applied to the commercial aspect. For example, Ma et al. [3] combined traditional taxi service with ridesharing to improve travel efficiency so that businesses could improve their income. Caulfield [4] revealed that both the society and environment benefits of ridesharing in terms of reductions in emission and energy. However, few studies integrated ridesharing into disaster response.

In recent years, people have noticed that ridesharing could improve disaster evacuation, and relevant examples have also emerged in the United States. For example, Lyft offered free trips to evacuation centers in Los Angeles and Ventura in 2017 [5]. Therefore, researchers have realized the importance of integrating ridesharing into disaster response, as advantages such as saving cost for individuals, providing flexibility for businesses, and optimizing transportation modes for society can improve disaster response efficiency. Currently, Wong and Shaheen [6] are the only researchers who conducted a relevant study. However, they merely distributed surveys to individuals to analyze the attitude to share resources during disasters. So far, no one has studied this topic systematically and provided any mathematical models. Therefore, there are still some gaps between ridesharing and disaster response since researchers focused on this topic ideologically rather than analyzed it theoretically.

Based on the motivation of improving disaster response efficiency, this paper aims to fill the above-mentioned gaps by developing a mathematical model. Contributions of this paper are threefold. First, we integrate ridesharing into the traditional evacuation model by developing a mixed-integer programming model for the ridesharing evacuation problem. Second, a real-world case study based on Houston is used to validate the feasibility of the proposed model. Finally, we explore the disaster evacuation efficiency in a series of instances, such as different ratios of vehicles participating in the model and the maximum evacuation time.

This paper is organized as follows. Section 2 presents relevant research about disaster evacuation and ridesharing. Section 3 proposes the mathematical formulation of ridesharing evacuation. Section 4 tests the validity of the ridesharing evacuation model and analyzes the efficiency of disaster response in multiple instances. Finally, some conclusions are given in Section 5.

## 2. Literature Review

The research literature review is divided into two categories: disaster evacuation and ridesharing, and because of the page limit, we only focus on the most relevant literature here.

### 2.1. Disaster Evacuation

Previous studies related to disaster evacuation have been conducted extensively. Many existing studies focused on public transport evacuation, such as moving carless people to shelters by public transport and exploring public transport plans in multiple disaster scenarios. For example, Sayyady and Eksioglu [1] proposed a mixed-integer linear program to determine the optimal evacuation route for individuals depending on public transport. Gao et al. [2] developed a robust two-stage transit-based evacuation model for large-scale disaster response. However, these studies assumed some meeting points before the evacuation without considering how to transfer individuals to these places. Besides, most of these studies focused on designing evacuation routes or allocating vehicles based on public transport rather than private vehicles. Meanwhile, methodologies such as statistical and probabilistic models, queuing theory, simulation, decision theory, and optimization models are proposed by previous researches. The optimization model is the most commonly used methodology because it can accurately represent each requirement in the highly complex evacuation process [7]. For example, Yi and Kumar [8] developed an optimization model for solving the logistics problem in disaster relief activities. Iakovou et al. [9] developed a linear programming model for emergency equipment logistics. Therefore, this paper aims to propose a mixed-integer programming model that individuals with vehicles can pick up others along their route to evacuation for disaster response.

### 2.2. Ridesharing

Ridesharing is an emerging topic and many researchers have studied it over the past decade. For example, Chan and Shaheen [10] presented the growth of ridesharing in North America. Agatz et al. [11] defined optimization problems in ridesharing and provided directions for future research by presenting the mechanism of ridesharing. Furuhata [12] revealed issues with the state-of-the-art ridesharing system and demonstrated the future direction of the ridesharing system. From these studies, we find that the most common application of previous studies focused on commercial aspects, such as optimizing transport methods to reduce congestion, decrease emission, and save resources. Very few researchers recently notice these advantages of ridesharing could also enhance disaster response. For instance, Wong and Shaheen [6] distributed surveys to individuals impacted by three major wildfires in California to analyze the attitude to share resources during disasters. However, they only investigated the feasibility rather than proposing any models to implement the integration. Therefore, there is still a gap in applying ridesharing to disaster response.

In conclusion, our study aims to propose a mixed-integer programming that optimizes the individuals to be evacuated quickly and people without vehicles could reach the public transport gathering places with the help of a ridesharing evacuation model. At the same time, we seek to improve resource utilization during disaster evacuation by exploring the performance of this model in different instances.

## 3. Model Formulation

In this paper, the idea of ridesharing is introduced into the traditional evacuation process to improve the efficiency of disaster response. Generally, participants involved in the ridesharing evacuation plan can be divided into two groups if they have their own vehicles to evacuate or not. When a disaster occurs, carless individuals need a ride for evacuation, while individuals who have vehicles can choose either evacuate to a gathering location directly for public transport or provide a ride to carless individuals along the way. Definitions of sets, parameters, and decision variables are given below.

**Sets**

| | |
|---|---|
| *R* | Locations of individuals who have vehicles |
| *H* | Locations of carless individuals |
| $S$ | Locations of gathering places where relief agencies (e.g., American Red Cross and FEMA) and local governments will provide the public transport for evacuation |

**Parameters**

| | |
|---|---|
| $d_i$ | Demand of individuals who need to be evacuated at location $i$ |
| $t_P$ | Time of each person to be loaded into a vehicle |
| $t_{ij}$ | Travel time from location $i$ to location $j$ |
| $t_{max}$ | Maximum time of individuals can be evacuated to a gathering place before disasters happen |
| $c^k$ | Capacity of vehicle $k$ |
| $M$ | Large positive number |

**Decision variables**

| | |
|---|---|
| $x_{ij}^k$ | If vehicle k travels from $i$ to $j$, this variable is 1, otherwise is 0 |
| $z_i$ | If any individual be evacuated at location $i$, this variable is 1, otherwise is 0 |
| $y_i^k$ | Number of individuals to be evacuated by vehicle $k$ at location $i$ |
| $u_i^k$ | Number of individuals to be evacuated by vehicle $k$ when leaves location $i$ |
| $v_i$ | Time at which vehicle leaves location $i$ |

The objective function is to maximize the number of individuals to be evacuated to a gathering place for public transport evacuation within a maximum time $t_{max}$. The mathematical model is formulated as follows.

$$Max \sum_{i \in R \cup H} \sum_{k \in R} y_i^k \qquad (1)$$

$$\sum_{j \in R \cup H \cup S} x_{ij}^i = z_i \quad \forall i \in R \qquad (2)$$

$$\sum_{i \in R \cup H} x_{ij}^k - \sum_{i \in H \cup S} x_{ji}^k = 0 \quad \forall j \in H, \forall k \in R \qquad (3)$$

$$\sum_{i \in R \cup H} \sum_{j \in S} \sum_{k \in R} x_{ij}^k = \sum_{i \in R} z_i \qquad (4)$$

$$y_i^k = d_i * z_i \quad \forall i \in R, \forall k \in R, If\ i = k \qquad (5)$$

$$\sum_{k \in R} y_i^k \geq z_i \quad \forall i \in H \qquad (6)$$

$$\sum_{k \in R} y_i^k \leq d_i * z_i \quad \forall i \in H \qquad (7)$$

$$y_j^k \leq M * \sum_{i \in R \cup H} x_{ij}^k \quad \forall j \in H, \forall k \in R \qquad (8)$$

$$u_i^k = y_i^k \quad \forall i \in R, \forall k \in R \qquad (9)$$

$$u_j^k \leq M * \sum_{i \in R \cup H} x_{ij}^k \quad \forall j \in H, \forall k \in R \qquad (10)$$

$$u_j^k \geq u_i^k + y_j^k - M * (1 - x_{ij}^k) \quad \forall i \in R \cup H, \forall j \in H, \forall k \in R \qquad (11)$$

$$v_j \geq v_i + t_{ij} + t_P * y_i^k - M * (1 - x_{ij}^k) \quad \forall i \in R \cup H, \forall j \in R \cup H, \forall k \in R \qquad (12)$$

$$v_i \leq t_{max} - t_{ij} * x_{ij}^k \quad \forall i \in R \cup H, \forall j \in S, \forall k \in R \qquad (13)$$

$$u_i^k \leq c^k \quad \forall i \in R \cup H, \forall k \in R \qquad (14)$$

$$x_{ij}^k = 0 \quad \forall i \in R, \forall j \in R \cup H \cup S, \forall k \in R, If\ i \neq k \qquad (15)$$

$$x_{ij}^k = 0 \quad \forall i \in R \cup H \cup S, \forall j \in R \cup H \cup S, \forall k \in R, If\ i = j \qquad (16)$$

$$u_i^k = 0 \quad \forall i \in R, \forall k \in R, If\ i \neq j \qquad (17)$$

$$y_i^k = 0 \quad \forall i \in R, \forall k \in R, If\ i \neq j \qquad (18)$$

$$x_{ij}^k \in \{0,1\},\ y_i^k \in \mathbb{Z}_+, z_i \in \{0,1\}, u_i^k \geq 0, v_i \geq 0, \quad \forall i,j \in R \cup H \cup S, \forall k \in R \qquad (19)$$

Constraint (2) limits that only vehicle $i$ can leave location $i$. Constraint (3) makes sure that the number of vehicles arriving at location $i$ is the same as the number of vehicles departing from that location. Constraint (4) indicates that every vehicle $k$ eventually arrives at a gathering place. Constraint (5) updates the number of individuals evacuated by vehicle $k$ at location $i$. Constraints (6) – (8) ensure that $y_i^k$ is updated when individuals need to be evacuated at location $i$. Constraint (9) updates the number of individuals each vehicle carries at the starting locations (locations in set $R$). Constraint (10) ensures that $u_j^k$ is updated only when vehicles travel from location $i$ to location $j$. Constraint (11) updates the capacity in each route. Constraint (12) updates the evacuation time for each route. Constraint (13) indicates that each vehicle $k$ must arrive at a gathering place before the deadline $t_{max}$. Constraint (14) indicates that every vehicle $k$ has a maximum capacity $c^k$. Constraints (15) – (18) indicate each vehicle $k$ must start from its starting location $i$. Finally, constraint (19) defines all decision variables.

## 4. Case Study
### 4.1. Information and dataset
A report in Natural Hazard Center [13] indicated that Houston is vulnerable to natural disasters like hurricane and flooding because the landscape of Houston is incredibly flat, and the geolocation of Houston is near a coast along with a warming ocean. Besides, people in such a big city like Houston do not need cars since there is a convenient public transport system. Therefore, the real-world case in Houston is used to test the ridesharing evacuation model. For comparison purposes, five different instances with multiple maximum evacuation times are discussed.

The specific case study area of this paper is the first ring of the Houston city, surrounded by Interstate 610 (I-610) freeway because Houston is a radial city from the perspective of urban planning where the inner area has a higher population density. According to the data provided by the Houston government and ArcGIS database[1], 14 households are randomly selected in this study area for experiments. Based on the information released by the Houston government[2], 8 safe gathering places are picked, as shown in Figure 1. Two examples of results are also shown in Figure 2 (a, b), detail explanations are given in the next section.

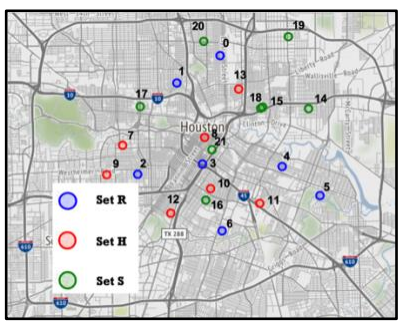
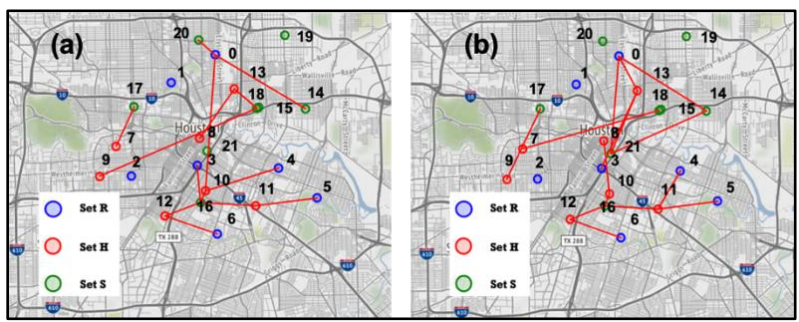

**Figure 1:** Case study area                **Figure 2:** Examples of results

Previous research [6] indicated that fewer than 50% of individuals were willing to share resources with others in a disaster. However, results may differ with different percentages of individuals who provide rides in our case. Therefore, we are interested in exploring multiple instances, which have different ratios of individuals in set $R$ and set $H$. To be more specific, we design a total of five instances: individuals in set $R$ account for 30%, 40%, 50%, 60% and 70% of the total participants respectively.

A set of parameters for the proposed ridesharing evacuation model is designed as follows. First of all, on the basis of Houston demographic data[3], the average household size of Houston residents is 2.76, so the number of individuals to be evacuated from each location is assumed to be 3. It is assumed that half of the vehicles in this case study have a capacity of 5-seats and others have a capacity of 7-seats. In order to make the case study more realistic, travel distance is calculated from Google Maps instead of Euclidean distance. Each individual is assumed to take 1 minute to get on the vehicle. According to the travel distance between places, the maximum evacuation time is set to 5, 7, 9, 11, 13, and 15 minutes respectively.

### 4.2. Results
The case study was carried out on a Dell Tower 3420 desktop with Intel (R) Core (TM) i7-7700 CPU @ 3.60GHz under a Windows 10 environment and Gurobi Optimizer 8.1. For all the instances, computational time ranges from 0.01 seconds to 251.94 seconds. Two examples of results are shown in Figure 2 (a, b) above, which demonstrate the optimal routes for the ridesharing evacuation plan when the maximum evacuation times are 13 and 15 minutes respectively. We find that the ridesharing evacuation model can solve the optimal evacuation route for different instances. In addition, we find that it is not necessary to increase the vehicles of disaster response because sometimes the optimal solution can be obtained with fewer vehicles. More detailed results are shown in our Mendeley database[4].

---

[1] ArcGIS dataset: http://mycity.houstontx.gov/home/
[2] Locations' information: https://www.naacp.org/campaigns/hurricane-harvey-2/
[3] Demographic data: https://datausa.io/profile/geo/houston-tx

Case study results show that the ridesharing evacuation model can evacuate all individuals in a limited time. For example, within the capacity of each vehicle, the ridesharing evacuation model can transfer all participants within 15 minutes in all instances, providing relief agencies plenty of preparation time for public transport evacuation. As the evacuation time and the ratio of individuals participating in the evacuation increases, more and more people will be successfully evacuated. The Department of Homeland Security's data[5] indicates the government will issue evacuation orders between 12 and 22 hours before a hurricane strike. Therefore, the ridesharing evacuation plan not only can ensure the people with vehicles evacuate safely but also help carless people to arrive at safe gathering places for public transport evacuation on time.

In this paper, the ridesharing evacuation plan also seeks to improve the efficiency of disaster evacuation. Therefore, two indicators are used for further analysis. One is evacuation percentage (EP), which indicates a percentage of the evacuated individuals in the total number of participants who need to be evacuated. The other is average travel distance (ATD), which represents the average travel distance per vehicle. Detailed results are shown in Figure 3 below.

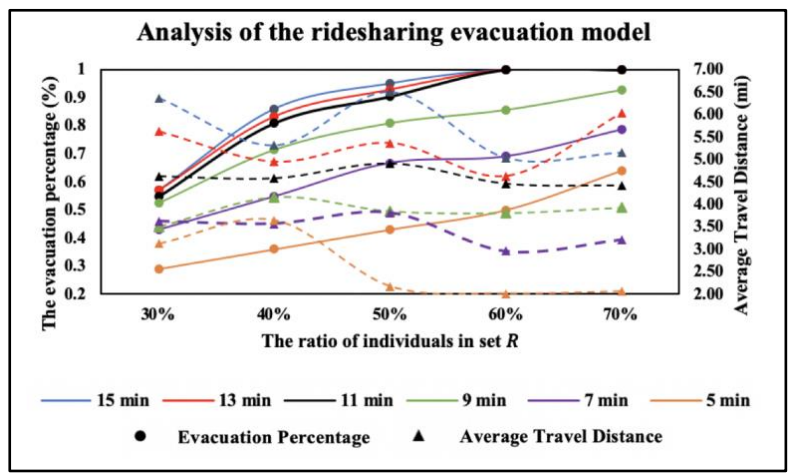

**Figure 3:** Analysis of the ridesharing evacuation model

Firstly, we analyze evacuation percentage (EP) in different instances. From Figure 3, we can see as the ratio of individuals in set *R* and the evacuation time increases, EP goes up. But the increasing rate of EP has changed over various evacuation times. For example, when the evacuation time is 5 minutes, as the ratio of individuals in set *R* increases, EP increases quickly. However, when the evacuation time is longer, i.e. more than 10 minutes, the increment of EP is flatter than the previous instance. Therefore, we find that it does not mean that the more individuals provide a ride, the better disaster evacuation is. In fact, we can obtain the optimal EP result though there is a lower ratio of individuals in set *R*. For example, the value of EP can achieve 100% when the ratio of individuals in set *R* is 60% within 13 or 15 minutes evacuation time.

Regarding the indicator ATD, from Figure 3 we can see that as the evacuation time increases, ATD grows. But we find that increasing the ratio of vehicles may not necessarily improve evacuation efficiency. For example, when the evacuation time is 5 minutes and the ratio is 30%, ATD is only 3.125 miles. When evacuation time does not change but the ratio is 40%, ATD increases to 3.57 miles. However, the ATD suddenly drops when ratio is 50% and then the results of ATD are very close when the ratio is 60% and 70% respectively.

In conclusion, the ridesharing evacuation model proposed in this paper performs well to evacuate individuals in the case study since within the capacity of each vehicle, the ridesharing evacuation model can evacuate all participants within 15 minutes in all instances. Besides, implementing ridesharing in disaster evacuation can improve the efficiency of disaster response as well as resource utilization because people without cars can be transferred to gathering places by people with cars. Meanwhile, the ridesharing evacuation model also demonstrates that it is not necessary to increase the resources of disaster response because the optimal solution can be obtained with fewer resources. Also, using too many vehicles may cause a problem like congestion and hinder the implementation of the evacuation. For instance, the case study area is the central area of Houston city, where the population and vehicle density are much higher than in other areas. Therefore, adding vehicles will cause more spaces to be occupied, such as road networks and parking spaces. Consequently, the ridesharing evacuation plan could conflict with public transport and result in negative impacts such as congestions of traffic and the suffering of individuals. Therefore, relief agencies should choose a specific evacuation strategy based on different disaster scales by this ridesharing evacuation model.

---

[4] Mendeley database: https://data.mendeley.com/datasets/7rh2ypk8yv/draft?a=32a7ef57-8cb6-4d1d-939b-92f8e4480404

[5] Homeland Security's data: https://www.ready.gov/hurricanes

## 5. Conclusion

This paper develops a mixed-integer programming model for maximizing the number of individuals to be evacuated to a gathering place within a maximum time $t_{max}$. In addition, to improve the efficiency of resource utilization, indicators evacuation percentage (EP) and average travel distance (ATD) are analyzed with different instances under a real-world case study in Houston. Results show that the ridesharing evacuation plan can evacuate all individuals in a limited time, and it is not necessary to increase the number of resources for disaster response because relief agencies should design a specific response strategy based on the disaster scales according to this model. It will not only provide disaster relief agencies and individuals with an efficient evacuation plan but also offer relief agencies insights in resource utilization.

The model proposed in this study belongs to the NP-hard problem. With the expansion of the case size, feasible solutions of the NP-hard problem will increase accordingly, thus posing a considerable challenge to the computation of the optimal solution. Since the model in our case is based on small sample size, future studies will focus on collecting more data to test multiple disasters to provide more practical conclusions.

## References

[1] F. Sayyady and S. D. Eksioglu, "Optimizing the use of public transit system during no-notice evacuation of urban areas," *Comput. Ind. Eng.*, vol. 59, no. 4, pp. 488–495, 2010.
[2] X. Gao, M. K. Nayeem, and I. M. Hezam, "A robust two-stage transit-based evacuation model for large-scale disaster response," *Meas. J. Int. Meas. Confed.*, vol. 145, pp. 713–723, 2019.
[3] S. Ma, Y. Zheng, and O. Wolfson, "T-share: A large-scale dynamic taxi ridesharing service," *Proc. - Int. Conf. Data Eng.*, pp. 410–421, 2013.
[4] B. Caulfield, "Estimating the environmental benefits of ride-sharing: A case study of Dublin," *Transp. Res. Part D Transp. Environ.*, vol. 14, no. 7, pp. 527–531, 2009.
[5] CBS NEWS, "Lyft offers free rides to evacuees of Southern California fires", *CBS Interactive Inc*. December 2017. Available: https://www.cbsnews.com/news/lyft-offers-free-rides-to-evacuees-of-southern-california-fires. [Accessed January 19, 2020].
[6] S. Wong, S. Shaheen, "Current State of the Sharing Economy and Evacuations: Lessons from California", *UC Office of the President: ITS reports*, 2019, pp. 1–50.
[7] A. M. Caunhye, X. Nie, and S. Pokharel, "Optimization models in emergency logistics: A literature review," *Socio-econ. Plann. Sci.*, vol. 46, no. 1, pp. 4–13, 2012.
[8] W. Yi and A. Kumar, "Ant colony optimization for disaster relief operations," *Transp. Res. Part E Logist. Transp. Rev.*, vol. 43, no. 6, pp. 660–672, 2007.
[9] E. Iakovou, C. M. Ip, C. Douligeris, and A. Korde, "Optimal location and capacity of emergency cleanup equipment for oil spill response," *Eur. J. Oper. Res.*, vol. 96, no. 1, pp. 72–80, 1997.
[10] N. Chan and S. Shaheen, "Ridesharing in North America: Past, Present, and Future," *Transp. Rev.*, vol. 32, no. 1, pp. 93–112, 2012.
[11] N. Agatz, A. Erera, M. Savelsbergh, and X. Wang, "Optimization for dynamic ride-sharing: A review," *Eur. J. Oper. Res.*, vol. 223, no. 2, pp. 295–303, 2012.
[12] M. Furuhata, M. Dessouky, F. Ordóñez, M. E. Brunet, X. Wang, and S. Koenig, "Ridesharing: The state-of-the-art and future directions," *Transp. Res. Part B Methodol.*, vol. 57, pp. 28–46, 2013.
[13] P. Berke, "Why is Houston so Vulnerable to Flooding?", *Natural Hazard Center*. September 2017. Available: https://hazards.colorado.edu/news/research-counts/part-i-why-is-houston-so-vulnerable-to-flooding. [Accessed January 19, 2020].